\documentclass[11pt,a4paper]{article}
\usepackage{graphicx}
\usepackage{amsfonts}
\usepackage{a4,amsmath,amssymb,here,epsf,amsfonts} \usepackage{a4,amsmath,amssymb,here,epsf,amsfonts}

\newtheorem{definition}{Definition}[section]
\newtheorem{theorem}{Theorem}[section]

\newtheorem{corollary}{Corollary}[section]

\newtheorem{lemma}{Lemma}[section]
\newtheorem{remark}{Remark}[section]
\newtheorem{conclusion}{Conclusion}[section]

\newenvironment{demo}{\endsloppypar\noindent\bgroup\small
            {\bf{Proof.} \indent}
            }{\samepage\null\hfill \mbox{$\rule{2mm}{2mm}$}
\quad\endsloppypar\egroup}
\begin{document}

\title{Euclidean Jordan algebras and generalized Krein parameters of  a strongly regular   graph}
\author{ Lu\'{\i}s Ant\'onio de Almeida Vieira\\
Faculdade de Engenharia\\
Universidade do Porto, CMUP\\
R. Dr Roberto Frias\\
4200-465 Porto, Portugal
}
\date{ }
\maketitle
\begin{abstract}
Let $\tau$ be a strongly $(n,p;a,c)$ regular   graph, such that $0<c<p<n-1,$ $A$  his matrix of adjacency and let ${\cal V}_{n}$ be the Euclidean space spanned by the powers of $A$ over the reals where the scalar product $\bullet|\bullet$ is defined by $x|y=\mbox{trace}(x \cdot y).$ In this work ones proves that ${\cal V}_{n}$ is an Euclidean Jordan algebra of rank 3 when one introduces in ${\cal V}_{n}$ the usual product of matrices.  Working inside the Euclidean Jordan algebra ${\cal V}_{n}$ and making use of the properties of $|A|^x$ one defines the generalized Krein parameters of the strongly $(n,p;a,c)$ regular graph $\tau$ and finally one presents necessary conditions over the parameters and the spectra of the $\tau$ strongly $(n,p;a,c)$ regular graph.
\end{abstract}
\section{Introduction}\label{rute0}
Euclidean Jordan algebras are more and more used in various branches of Mathematics.  For example, one may cite the application of this theory to  statistics \cite{massam98}, to interior point methods\cite{Faybusovich97a, Faybusovich97b,guler96}, combinatorics \cite{bannai_ito84} \cite{dcardosovieira}\cite{pech98}.
Among several readable texts about or which includes
Euclidean Jordan algebras, we may refer \cite{faraut94} and \cite{koecher99}. Previous works
on strongly regular graphs included in Euclidean Jordan algebras were developed aiming to obtain the
minimal coherent algebra which includes a particular set of matrices and then its symmetric homogeneous
coherent $3$-dimension subalgebras (all of them defining strongly regular graphs, as it is well known)
\cite{pech98}.

This report is organized as follows.In the next section one presents the basic results of Euclidean Jordan algebras. In section \ref{sec0} on presents a brief introduction on strongly regular graphs.
 In section \ref{sec2} one associates an Euclidean Jordan Algebra ${\cal V}$ of rank three to the matrix of adjacency of a strongly regular  graph.
In section \ref{sec3} one defines the power of $|A|^x$ inside the Euclidean Jordan algebra ${\cal V}$ associated to a strongly regular  graph $\tau$ and one obtains his coordinates in the basis $\{I_{n},A,E_{1}\}$ where $A$  is the matrix of adjacency of $\tau$ and $E_{1}$ is a primitive idempotent of the unique Jordan frame of the Euclidean Jordan algebra ${\cal V}$ associated to $A.$
In section \ref{clara}  one defines the generalized Krein parameters of a strongly regular graph.  Finally  one establishes necessary conditions for the existence of a strongly $(n,p;a,c)$  regular graph.
\section{Basic results on Euclidean Jordan algebras}\label{sec1}
Consider a real finite and $n-$dimensional vectorial space $V$ with a multiplication $\circ$ such that the map
$(x,y) \rightarrow x \circ y$ is bilinear. If $\forall x \in V \;\; (x \circ x) \circ x = x \circ (x \circ x)$
one says that $V$ is a power associative algebra. In a power associative algebra, $V$, for all $x \in V$ and
for all $p,q \in \mathbb{N}$ $x^p \circ x^q=x^{p+q}.$ If $\textbf{e} \in V$ is such that
$\forall x \in V \;\; x \circ \textbf{e}=\textbf{e} \circ x = x,$ then $\textbf{e}$ is called the unit of $V$.\\
Let $V$ be a power associative algebra with unit element $\textbf{e}$. For all $x \in V$ let $k$ be the least
natural number such that $\{e,x,x^2,\ldots,x^{k}\}$ is linear dependent. Then $k$ is the rank of $x$ and one
writes $rank(x)=k$. One defines the rank of $V$ as being the natural number
$r=\mbox{rank}(V)=\mbox{max}\{rank(x):x\in V\}.$
A element $x \in V$ is regular if its rank is equal to the rank of the algebra. Given a regular element, $x$,
in a power associative algebra, $V$, with unit element, $\textbf{e}$, and rank $r$, since the set
$\{\textbf{e},x,x^2,\ldots,x^r\}$ is  linear dependent and the set $\{\textbf{e},x,x^2,\ldots,x^{r-1}\}$ is linear
independent, one may conclude that there are $r$ unique real numbers, $a_1(x),$ $a_2(x),$ $\ldots,$ $a_r(x)$, such
that
\begin{eqnarray}
x^r-a_1(x)x^{r-1}+\ldots+(-1)^ra_r(x)\textbf{e} = 0,\label{marca222222.1}
\end{eqnarray}
where $0$ is the null vector of $V$. Taking in account (\ref{marca222222.1}), the polynomial
\begin{equation}
p(x,\lambda) = \lambda^r-a_1(x)\lambda^{r-1}+\ldots+(-1)^ra_r(x) \label{marca_21}
\end{equation}
is called the characteristic polynomial of $x$, where the coefficients $a_i(x)$ are polynomial functions in the
coordinates of $x$ in a fixed basis of $V$. The definition of the characteristic polynomial may be extended to any
element of $V$. Indeed, since the set of regular elements of $V$ is a dense set in $V$ \cite{faraut94}, if  $x \in V$
then there exists a sequence  $\{ x_n \}_{n \in \mathbb{N}}$ of regular elements of $V$ converging to $x$. Defining
$\lim_{n\rightarrow \infty}a_i(x_n) = a_i(\lim_{n\rightarrow \infty}x_n)=a_i(x),$ one obtains the characteristic
polynomial of a non regular element as being the polynomial (\ref{marca_21}).\\
Let $x$ be a regular element of $V$. If one considers the real subalgebra of $V$, $\mathbb{R}[x]$, spanned by
$\{e,x, \ldots, x^{r-1}\}$, then the restriction of the application $L(x)$, such that $L(x)y=x \circ y$, to
$\mathbb{R}[x]$ denoted by $L_{0}(x)$ has the matrix representation, in the basis $B = \{ e,x, \ldots, x^{r-1} \}$,
given by
\begin{eqnarray*}
        M_{L_0(x)} & = & \left[\begin{array}{ccccc}
                                            0   &    0   & \cdots &    0   & (-1)^{r-1}a_r(x)\\
                                            1   &    0   & \cdots &    0   & (-1)^{r-2}a_{r-1}(x)\\
                                            0   &    1   & \cdots &    0   & (-1)^{r-3}a_{r-2}(x)\\
                                         \vdots & \vdots & \ddots & \vdots & \vdots\\
                                            0   &    0   & \cdots &    1   & a_1(x)
                              \end{array}
                          \right],
\end{eqnarray*}
whose characteristic polynomial is the polynomial (\ref{marca_21}).\\
From now on, $p(x,\lambda)$ stands for the characteristic polynomial of an element $x$ of a power associative algebra.
Easily one can see that when $x$ is a  regular element, $a_r(x)$ is equal to determinant of the matrix $M_{L_0(x)}$ and
$a_1(x)$ is equal to the trace of the matrix $M_{L_0(x)}$. For this reason, in a real power associative algebra $, V,$
with unit element, one defines for each element $x$ of $V$ the trace and the determinant of $x$ and one denotes them,
respectively, by $\mbox{tr}(x)$ and $\det(x)$, as being, respectively, the coefficients $a_1(x)$ and $a_r(x)$ of the
polynomial (\ref{marca_21}). As an example of real power associative algebras with unit element one may refer the real
Jordan algebras with unit element which one describes next.

\begin{definition}
Let  $V$ be a real finite dimensional vector, with the operation of multiplication of vectors, $\circ$, determined by
the bilinear application $(x,y)\rightarrow x\circ y$. One says that $V$ is a Jordan algebra if
\begin{itemize}
\item[i)]{$x \circ y=y\circ x$},
\item[ii)]{$x\circ (x^2\circ y)=x^2\circ (x\circ y),$ where $x^2=x\circ x$}.
\end{itemize}
\end{definition}

From now on, when one refers to a Jordan algebra, $V$, one admits that $V$ is a real finite dimensional algebra
and with unit element $\textbf{e}$. Since Koecher \cite{koecher99} proves that if $V$ is a Jordan algebra then
$V$ is a power associative algebra, then to each element  $x \in V$ one may associate the respective characteristic
polynomial. Thus, if $r$ is the characteristic of $V$ and the characteristic polynomial of $x$ is the polynomial
(\ref{marca_21}), then one defines the determinant and the trace of $x$, which are denoted, respectively, $det(x)$
and $tr(x)$, by the equalities $\mbox{tr}(x) = a_1(x)$ and $\det(x) = a_r(x)$. Given a Jordan algebra one says that
the element $x$ is invertible if $\exists y \in \mathbb{R}[x]$
\footnote{Subalgebra of $V$ generated by $\textbf{e}$ and $x$.}
such that $x \circ y = \textbf{e}$, and then $y$ is designated the inverse of $x$ and denoted by $x^{-1}$.


A real Jordan algebra $V$ is Euclidean if there is a scalar product $<\bullet,\bullet>$ such that $<u \circ v,w>=<v,u \circ w>.$
Additionally, two elements $c,d \in V$ are orthogonal, relatively to the algebra $V$, if $c \circ d=0.$ On the other
hand, assuming that $\textbf{e}$ is the unit element of $V$, $c \in V$ is an idempotent if $c^{2}=c$, and one says
that $\{c_1,c_2,\ldots,c_k\}$ is a complete system of orthogonal idempotents if
\begin{eqnarray*}
     (i)   &\! & c^{2}_i=c_i \;\; \forall i \in \{1, \ldots , k \}, \\\!
      (ii) &\! &  c_i \circ c_j=0 \;\; \forall i \not= j,\\\!
     (iii) &\! & c_1 + c_2+ \cdots + c_k = \textbf{e}.
\end{eqnarray*}
An idempotent $c$ is primitive if it is not the sum of two other idempotents. One says that $\{c_1,c_2,\ldots,c_k\}$
is a complete system of orthogonal primitive idempotents or a Jordan frame, if $\{c_1,c_2,\ldots,c_k\}$ is a complete
system of orthogonal idempotents such that each idempotent is primitive.

\begin{theorem}(\cite{faraut94},p.43)\label{teorema_242}
Let $V$ be a Euclidean Jordan algebra with unit element $\textbf{e}.$ For $x\in V$ there exist unique real numbers $k$ unique
real numbers, $\lambda_1,\lambda_2,\ldots,\lambda_k$, all distinct, and a unique complete system of orthogonal
idempotents, $\{ c_1,c_2,\ldots,c_k \}$, such that
\begin{eqnarray}
x=\lambda_1c_1+\lambda_2c_2+\cdots+\lambda_kc_k. \label{marca_28}
\end{eqnarray}
Additionally, $c_j \in \mathbb{R}[x]$, for  $j=1,\ldots, k$.
\end{theorem}

The numbers $\lambda_j$ of (\ref{marca_28}) are called the eigenvalues of $x$ and $x = \sum^k_{i=1} {\lambda_i}c_i$
is the spectral decomposition of $x$.\\
If $V$ is a Euclidean Jordan algebra with unit element $\textbf{e}$ and characteristic $r$, and $\textbf{c}$ is a
primitive idempotent of $V$, then $\mbox{tr}(c)=1$ and, therefore, one may conclude that $\mbox{tr}(\textbf{e})=r.$

\begin{theorem}(\cite{faraut94},p.44)\label{teorema_243}
Let $V$ be an Euclidean Jordan algebra with characteristic $r$ and unit element $\textbf{e}$. Then, for each $x \in V,$
there is a Jordan frame  $\{c_1,c_2,\ldots,c_r\}$ and  real numbers $\lambda_1,\lambda_2,\ldots,\lambda_r$
such that
\begin{eqnarray*}
       x & = & \sum^r_{j=1}{\lambda_jc_j}.
\end{eqnarray*}
The  $\lambda`_i$s, together  with their multiplicities, are uniquely  determined by $x$. Additionally one verifies
that
\begin{eqnarray*}
  \mbox{tr}(x) &\! = & \sum^r_{j=1}{\lambda_j}.
\end{eqnarray*}

\end{theorem}
\section{Brief introduction to strongly regular graphs}\label{sec0}
A graph $G$, is a pair of sets $(V(G),E(G)),$ where $V(G)$ denotes the nonempty set of vertices and
$E(G)$ the set of edges. An element of  $E(G)$, which the endpoints are the vertices $i$ and $j$, is
denoted by $ij$ and, in such case, we say that the vertex $i$ is adjacent to the vertex $j.$
If $v \in V(G),$ then we call neighborhood of $v$ the vertex set denoted by
$N_{G}(v) = \{w: vw \in E(G)\}.$  The complement of $G,$ denoted by $\bar{G}$, is such that
$V(\bar{G})=V(G)$ and $E(\bar{G}) = \{ij: i,j \in V(G) \wedge ij \not \in E(G)\}.$ A graph of order
$n$ in which all pairs of vertices are adjacent is called a complete graph and denoted by $K_{n}.$
On the other hand, when there are no pair of adjacent vertices the graph is called a null graph and
denoted by $N_{n}$. In this paper we assume that $G$ is of order $n>1,$ i.e., $|V(G)|=n>1$ and we deal
only with simple graphs, that is, graphs with nor loops (edges with both ends in the same vertex) neither
parallel edges (more than one edge between the same pair of vertices). Consider the set of vertices of
the graph $G,$ $X = \{x_{1}, x_{2}, \ldots, x_{k}\},$ and suppose that
$x_{i}x_{i+1} \in E(G) \; \forall i \in \{1, \ldots, k-1\}.$ If the vertices of $X$ are all
distinct then we say that $X$ induces a path (and if they are all distinct but $x_{1}$ and $x_{k}$, which
are equal, then we say that $X$ induces a cycle). If there is no $i \in \{1, \ldots, k-2\}$ and $j>1$
such that $x_{i}x_{i+j} \in E(G),$ then $X$ is denoted by $P_{k}$ ($C_{k}$) and we say that it is a
chordless path (cycle). The adjacency matrix of the graph $G$ is the matrix
$A_{G} = \left(a_{ij}\right) _{n\times n}$ such that
$$
a_{ij}=\left\{\begin{array}{cl}
                     1 & ,\mbox{if } ij \in E(G) \\
                     0 & ,\mbox{otherwise.}
              \end{array}
       \right.
$$
The number of neighbors of $v \in V(G)$ will be denoted by $d_{G}(v)$ and called degree
of $v.$ If $G$ is such that $\forall v \in V(G) \;\; d_{G}(v)=p$ then we say that $G$ is $p$-regular.
A graph $G$ is called strongly regular if it is regular, not complete, not null and, given any two
distinct vertices $i,j \in V(G)$, the number of vertices which are neighbors of both $i$ and $j$
depends on whether $i$ and $j$ are adjacent or not. When $G$ is a strongly regular graph of order $n$
which is $k$-regular and any pair of adjacent vertices have $p$ common neighbors and any two distinct
non-adjacent vertices have $q$ common neighbors, then we say that $G$ is a $(n,p;a,c)$-strongly regular
graph. The chordless cycle on five vertices $C_{5}$ is an example of a strongly regular graph which is
(5,2;0,1)-strongly regular.\\
It is well known (see for instance \cite{Godsil2001} that a graph $G$ which is not complete and
not null is strongly regular iff $A_{G}^{2}$ is a linear combination of $A_{G}$, $I_{n}$ and
$J_{n}$, that is, $\exists \tau_{1}, \tau_{2}, \tau_{3} \in \mathbb{R}$ such that
\begin{eqnarray}
A_{G}^{2}&  = & \tau_{1}I_{n} + \tau_{2}A_{G} + \tau_{3}J_{n} \label{regularidade_forte},
\end{eqnarray}
where $J_{n}$ denotes the all ones square matrix of order $n$. Therefore, since
$A_{\bar{G}} = J_{n}-I_{n}-A_{G}$ we may conclude that a graph is strongly regular if and only if
its complement is also strongly regular.From now on one will use $A$ for representing the matrix of adjacency of a strongly regular graph.

Let $\tau$ be a strongly $(n,p;a,c)$ regular graph such that $0<c<p<n-1,$ $A$ his matrix of adjacency and ${\cal V}_{n}$ be the Euclidean space over the reals of the linear combinations of the  powers of $A$  with exponent  in $\mathbb{N}_{0}$ where the sum of vectors is the usual sum of matrices, the product of a vector  by a scalar is the product of matrix by a real number and the scalar product is $\bullet|\bullet$ defined by $x|y=\mbox{trace}(x\cdot y).$
\section{Euclidean Jordan algebra associated to the matrix  of adjacency  of a strongly regular  graph}\label{sec2}
Since the powers of $A^n$ with $n$ a natural number commute among each other then $x \cdot y=y \cdot x,\: \forall x,y \in {\cal V}_{n}.$
For the same reason and by the  fact that the product of matrices his associative, one concludes that
\begin{itemize}
\item[i)]
$x\cdot(x^2 \cdot y)=x^2\cdot(x \cdot y),\:\forall x,y\in {\cal V}_{n},$
\item[ii)]
$<L(x)y,z>=<y,L(x)z>,\: \forall x,y,z\in {\cal V}_{n}.$
\end{itemize}
Thus, the structure $({\cal V}_{n},\cdot)$ where $\cdot$ is the usual product of matrices is an Euclidean Jordan algebra.

One now introduces some notation.
Let $r$ and $s$ be the real numbers defined by
\begin{eqnarray*}
r&\!=&\frac{a-c+\sqrt{(a-c)^2+4*(p-c)}}{2},\\\!
s&\!=&\frac{a-c-\sqrt{(a-c)^2+4*(p-c)}}{2}.
\end{eqnarray*}
Since the linear operator $L(A)$ is a symmetric linear operator from ${\cal V}_{n}$ to ${\cal V}_{n}$ with only three distinct eigenvalues, namely $p,r$ and $s$ one concludes that
$S=\{E_{1},E_{2},E_{3}\}$ is the unique complete system of orthogonal   idempotents of ${\cal V}_{n}$ associated to $A,$ where
\begin{eqnarray*}
E_{1}&\!=&\frac{A^2-(r+s)A+rsI_{n}}{(p-r)(p-s)}=\frac{J}{n},\\\!
E_{2}&\!=&\frac{A^2-(p+s)A+psI_{n}}{(r-s)(r-p)},\\\!
E_{3}&\!=&\frac{A^2-(p+r)A+prI_{n}}{(s-r)(s-p)}
\end{eqnarray*}
and $J_{n}$ is the matrix of ones.
In the following one shows that $S$ is a Jordan frame of the Euclidean Jordan algebra ${\cal V}_{n}$ that is a basis of ${\cal V}_{n}$ with the aim of proving that ${\cal V}_{n}$ is an Euclidean Jordan algebra with rank three.
Since the spectral decomposition of $A$ is
$A=pE_{1}+rE_{2}+sE_{3}$ one concludes that
$$A^k=p^kE_{1}+r^kE_{2}+s^kE_{3},\:\forall k\in \mathbb{N}_{0}$$ and so the set $S$  spans ${\cal V}_{n}$ and  since
 $\{E_{1},E_{2},E_{3}\}$ is a linear independent set of ${\cal V}_{n}$ then
$<E_{1},E_{2},E_{3}>$ is a basis of ${\cal V}_{n}.$
Therefore  one may affirm that the dimension of  ${\cal V}_{n}$ is three and that $S$ is a  Jordan frame of ${\cal V}_{n}.$
Now let us prove that $\mbox{rank}(A)=3.$ Since $A$ is the matrix of  adjacency  $\tau$ then
$$(A-pI_{n})(A^2-(s+r)A+srI_{n})=0.$$ So, after some calculations it follows that
$$A^3=A^2(s+r+p)-(sr+(s+r)p)A+psrI_{n}.$$ Thus, one may conclude that $\{I_{n},A,A^2,A^3\}$ is a linear dependent set.
We will now prove that $\{I_{n},A,A^2\}$ is a free set of ${\cal V}_{n}.$  Let $\alpha, \beta$ and $\gamma$ be real numbers. Since
$$\left|\begin{array}{lll}
1&p&p^2\\
1&r&r^2\\
1&s&s^2\\
\end{array}
\right|=(r-p)(s-p)(s-r)\not=0$$
and
\begin{eqnarray*}
\alpha I_{n}+\beta A+\gamma A^2&\!=&0\\\!
&\!\Updownarrow&\\\!
\alpha E_{1}+\alpha E_{2}+\alpha E_{3}+\beta(pE_{1}+rE_{2}+sE_{3})+\gamma (p^2E_{1}+r^2E_{2}+s^2E_{3})&\!=&0\\\!
&\!\Updownarrow&\\\!
(\alpha +p \beta+p^2\gamma)E_{1}+(\alpha+r \beta+r^2 \gamma)E_{2}+(\alpha+s \beta+s^2 \gamma)&\!=&0\\\!
&\!\Updownarrow&\\\!
\left\{\begin{array}{lll}
\alpha +p \beta+p^2\gamma&=&0\\
\alpha+r \beta+r^2 \gamma&=&0\\
\alpha+s \beta+s^2 \gamma&=&0\\
\end{array}
\right.\\\!
&\!\Updownarrow&\\\!
\alpha=0\wedge \beta=0\wedge \gamma&\!=&0
\end{eqnarray*}
then  $\{I_{n},A,A^2\}$ is a free set of  ${\cal V}_{n},$ and therefore  the least natural number $k$ such that  the set $\{I_{n},A,A^2,\cdots,A^k\}$ is linear dependent is three and so $\mbox{rank}(A)=3.$ Since the dimension of ${\cal V}_{n}$ is three it follows that $\mbox{rank}({\cal V}_{n})=3.$
\section{Properties of the Euclidean Jordan algebra associated to a strongly regular  graph}\label{sec3}
Let $\tau$ be a strongly regular  graph, ${\cal V}_{n}$ the Euclidean Jordan algebra associated to $\tau$ and $S=\{E_{1},E_{2},E_{3}\}$ the unique Jordan frame associated to the matrix of adjacency of $\tau.$
Making use of the spectral decomposition of $I_{n}$ and of $A$ one concludes that the set $\{I_{n},A,E_{1}\}$ is a linear independent set of ${\cal V}_{n}.$  Since ${\cal V}_{n}$ is 3-dimensional then one may conclude that $\{I_{n},A,E_{1}\}$ is a basis of ${\cal V}_{n}.$
In this section one considers the basis ${\cal B}=\{I_{n},A,E_{1}\}$ and we make use of the scalar product $<x,y>=\mbox{tr}_{{\cal V}_{n}}(x\cdot y),\:\forall x,y\in {\cal V}_{n}$ where now the notation $\mbox{tr}_{{\cal V}_{n}}(u)$ represents the trace  of the vector $u$ in the sence of the Euclidean Jordan algebra ${\cal V}_{n}.$ In first place, one defines the power of $|A|^x$ for any real number $x$ and next one  calculates  the coordinates of $|A|^x$ in the basis ${\cal B}.$
Since $S$ is an Jordan frame of ${\cal V}_{n}$ it follows that $\mbox{tr}_{{\cal V}_{n}}(E_{i})=1,\:\forall i=1,\cdots 3.$  Let $x$ be a real number.
One defines $|A|^{x}$ by the equality
$$|A|^{x}=p^{x}E_{1}+r^{x}E_{2}+|s|^{x}E_{3}.$$
 Now one determines the coordinates of $|A|^{x}$ in the basis ${\cal B}.$
So let $\alpha_{x}, \beta_{x}$ and $\gamma_{x}$ be real numbers
such that
\begin{eqnarray}
|A|^{x}&=&\alpha_{x}I_{n}+\beta_{x}A+\gamma_{x}E_{1}.\label{eq3}
\end{eqnarray}
Using the spectral decomposition of $|A|^x$ it follows that
\begin{eqnarray}
p^{x}E_{1}+r^{x}E_{2}+|s|^xE_{3}&=&\alpha_{x}I_{n}+\beta_{x}A+\gamma_{x}E_{1}\label{eq4}.
\end{eqnarray}
Now let us consider in ${\cal V}_{n}$ the scalar product $<\bullet,\bullet>$ defined by $$<x,y>=\mbox{tr}_{{\cal V}_{n}}(xy),\:\forall x,y\in {\cal V}_{n}.$$
Multiplying both members of (\ref{eq4}) successively by $E_{1},E_{2}$ and $E_{3}$ and applying $\mbox{tr}_{{\cal V}_{n}}$ to both members  of (\ref{eq4})
one obtains the system
\begin{eqnarray}
\left\{
\begin{array}{lll}
r^{x}&=&1\alpha_{x}+r\beta_{x}+0\gamma_{x}\\
|s|^{x}&=&1\alpha_{x}+s\beta_{x}+0\gamma_{x}\\
p^{x}&=&1\alpha_{x}+p\beta_{x}+	\gamma_{x}
\end{array}
\right..\label{sistema}
\end{eqnarray}
Since
\begin{eqnarray*}
\left|\begin{array}{lll}
1&r&0\\
1&s&0\\
1&p&1
\end{array}
\right|&=&s-r\not=0
\end{eqnarray*}
then one may conclude that the system (\ref{sistema}) is a Cramer  system.
Therefore solving (\ref{sistema}) one obtains:
$$
\begin{array}{lll}
\alpha_{x}=\frac{(p-c)(r^{x-1}+|s|^{x-1})}{r-s},&
\beta_{x}=-\frac{|s|^{x}-r^{x}}{r-s},&
\gamma_{x}
=p^x-r^x+(p-r)\frac{|s|^x-r^x}{r-s}.
\end{array}
$$
Finally, one may write
 \begin{eqnarray}
|A|^x&\!=&\frac{(p-c)(r^{x-1}+ |s|^{x-1})}{r-s}I_{n}-\frac{|s|^{x}-r^{x}}{r-s}A+(p^{x}-r^x+(p-r)\frac{|s|^{x}-r^{x}}{r-s})E_{1}.\nonumber\\\label{finalmente}
\end{eqnarray}
Now let suppose that  $r,|s|$ and $p$ are distinct numbers.
Since $$|A|^{x}=p^xE_{1}+|r|^xE_{2}+|s|^{x}E_{3}$$ then rewriting (\ref{finalmente}) one obtains
 \begin{eqnarray}
p^xE_{1}+r^xE_{2}+|s|^xE_{3}&\!=&\frac{(p-c)(r^{x-1}+ |s|^{x-1})}{r-s}I_{n}-\frac{|s|^{x}-r^{x}}{r-s}A\nonumber\\\!
&\!+&(p^{x}-r^x+(p-r)\frac{|s|^{x}-r^{x}}{r-s})E_{1}.\label{finalmente1}
\end{eqnarray}
Working in the real  vector space ${\cal F}$ of real functions with the usual operations of addition of functions and scalar multiplication of a function by a real  and equalizing the coefficients of $p^x, r^x$ and of $|s|^x$ in  both members of (\ref{finalmente1}) one obtains:
\begin{eqnarray*}
E_{1}&=&\frac{J}{n},\\
E_{2}&=&\frac{|s|}{r-s}I_{n}+\frac{1}{r-s}A+\frac{s-p}{r-s}\frac{J}{n},\\
E_{3}&=&\frac{r}{r-s}I_{n}-\frac{1}{r-s}A+\frac{p-r}{r-s}\frac{J}{n}.
\end{eqnarray*}
Using the basis $\{I_{n},A,J_{n}-A-I_{n}\}$ of ${\cal V}_{n}$ it follows that
\begin{eqnarray}
E_{1}&=&\frac{r-s}{n(r-s)}I_{n}+\frac{r-s}{n(r-s)}A+\frac{r-s}{n(r-s)}(J_{n}-A-I_{n}),\nonumber\\
E_{2}&=&\frac{|s|n+s-p}{n(r-s)}I_{n}+\frac{n+s-p}{n(r-s)}A+\frac{s-p}{n(r-s)}(J_{n}-A-I_{n}),\nonumber\\
E_{3}&=&\frac{rn+p-r}{n(r-s)}I_{n}+\frac{-n+p-r}{n(r-s)}A+\frac{p-r}{n(r-s)}(J_{n}-A-I_{n})\label{ver}.
\end{eqnarray}
One now suppose that   $r, s$ and $p$ are distinct real numbers and that $r=|s|.$
Let prove that the expressions for the idempotents $E_{1},E_{2}$ and $E_{3}$ are the same as those  presented on (\ref{ver}) when $p,$$r$ and $|s|$ are not distinct. Let $\epsilon$ be a  sufficient small real number such that the set $\{(p+\epsilon)^x,(r+\epsilon)^x,|s+\epsilon|^x\}$ is a linear independent set of the real vector space ${\cal F}.$ Let now work in the Euclidean Jordan Algebra generated  by the powers of the matrix $A+\epsilon I_{n}.$ Then proceeding in the same way as when one has deduced (\ref{finalmente1}) with $p,r$ and $|s|$  distinct one  obtains:
\begin{eqnarray}\label{alg999}
\begin{array}{l}
(p+\epsilon)^{x}E_{1}+(r+\epsilon)^{x}E_{2}+|s+\epsilon|^{x}E_{3}=\\
|s+\epsilon|(r+\epsilon)\frac{(r+\epsilon)^{x-1}+ |s+\epsilon |^{x-1}}{r-s}I_{n}
-\frac{|s+\epsilon |^{x}-(r+\epsilon )^{x}}{r-s}(A+\epsilon I)+\\
+((p+\epsilon)^{x}-(r+\epsilon)^x+(p-r)\frac{|s+\epsilon|^{x}-(r+\epsilon)^{x}}{r-s})E_{1}.
\end{array}
\end{eqnarray}
Now working in the real  vector space ${\cal F}$  and equalizing the coefficients of $(p+\epsilon)^x,(r+\epsilon)^x$ and $|s+\epsilon|^x$ in  both members of (\ref{alg999}) one obtains:
\begin{eqnarray*}
E_{1}&=&\frac{J}{n},\\
E_{2}&=&\frac{|s+\epsilon|}{r-s}I_{n}+\frac{1}{r-s}(A+\epsilon I_{n})+\frac{s-p}{(r-s)}\frac{J}{n},\\
E_{3}&=&\frac{r+\epsilon}{r-s}I_{n}-\frac{1}{r-s}(A+\epsilon I_{n})+\frac{p-r}{(r-s)}\frac{J}{n}.
\end{eqnarray*}
Using the basis $\{I_{n}, A, J_{n}-A-I_{n}\}$ of ${\cal V}_{n}$ it follows that
\begin{eqnarray*}
E_{1}&=&\frac{J}{n},\\
E_{2}&=&\frac{|s+\epsilon|n+s-p}{n(r-s)}I_{n}+\frac{n+s-p}{n(r-s)}A
+\frac{n \epsilon}{n(r-s)}I_{n}+\frac{s-p}{(n(r-s)}(J_{n}-A-I_{n}),\\
E_{3}&=&\frac{(r+\epsilon)n+p-r}{n(r-s)}I_{n}+\frac{-n+p-r}{n(r-s)}A-\frac{n \epsilon }{n(r-s)}I_{n}+\frac{p-r}{n(r-s)}(J_{n}-A-I_{n}).
\end{eqnarray*}
Then making $\epsilon $ converge  to zero one obtains:
\begin{eqnarray}
E_{1}&=&\frac{J}{n},\nonumber\\
E_{2}&=&\frac{|s|n+s-p}{n(r-s)}I_{n}+\frac{n+s-p}{n(r-s)}A+\frac{s-p}{n(r-s)}(J_{n}-A-I_{n}),\nonumber\\
E_{3}&=&\frac{rn+p-r}{n(r-s)}I_{n}+\frac{-n+p-r}{n(r-s)}A+\frac{p-r}{n(r-s)}(J_{n}-A-I_{n}).\label{marilia}
\end{eqnarray}
Now one will establish necessary conditions for the existence of a strongly $(n,p;a,c)$ regular graph $\tau$ working inside the Euclidean Jordan algebra ${\cal V}_{n}$ associated to $\tau$ and using the theorem \ref{interlace}.
\section{Generalized Krein parameters of a strongly regular graph}\label{clara}
In this section one  will define the generalized Krein parameters of the strongly regular graph $\tau,$   establish some properties of these parameters and finally one will deduce necessary conditions for the existence of a strongly $(n,p;a,c)$ regular graph. But first one presents the theorems \ref{interlace},\ref{kronecker11}, the lemma \ref{kronecker12} and some notation.
\begin{theorem}(\cite{jnvanlint},p.439)\label{interlace}
Let $A$ be a symmetric matrix of order $n$ with eigenvalues
$$\lambda_{1}\geq \lambda_{2}\geq \cdots \geq \lambda_{n}.$$
Suppose $N$ is an $m \times n$ real matrix such that $NN^T=I_{m},$ so $m\leq n.$ Let $B=NAN^T,$ and let $$\mu_{1}\geq \mu_{2}\geq \cdots\geq \mu_{m}$$ be the eigenvalues of $B.$ Then the eigenvalues of $B$ interlace those of $A,$ in the sense that
$$\lambda_{i}\geq \mu_{i}\ge \lambda_{n-m+i}$$
for $i=1,2,\cdots,m.$
\end{theorem}
From (\ref{marilia})  it follows
that
$$
\begin{array}{l}
E_{1}=\frac{(r-s)I_{n}+(r-s)A+(r-s)(J_{n}-A-I_{n})}{n(r-s)},\\
E_{2}
=\frac{|s|n+s-p}{n(r-s)}I_{n}
+\frac{n+s-p}{n(r-s)}A
+\frac{s-p}{n(r-s)}(J_{n}-A-I_{n}),\\
E_{3}
=\frac{rn+p-r}{n(r-s)}I_{n}
+\frac{-n+p-r}{n(r-s)}A+
\frac{p-r}{n(r-s)}(j_{n}-A-I_{n}),\\
E_{1}+E_{2}=\frac{|s|n+r-p}{n(r-s)}I_{n}+\frac{n+r-p}{n(r-s)}A+\frac{r-p}{n(r-s)}(J_{n}-A-I_{n}),\\
E_{1}+E_{3}=\frac{rn+p-s}{n(r-s)}I_{n}+\frac{-n+p-s}{n(r-s)}A+\frac{p-s}{n(r-s)}(J_{n}-A-I_{n}),\\
E_{2}+E_{3}=	\frac{n-1}{n}I_{n}-\frac{1}{n}A-\frac{1}{n}(J_{n}-A-I_{n}).
\end{array}
$$
Now one will introduce some notation.
\begin{definition}
{\rm
let $A$ and $B$ be matrices of ${\cal M}_{m\times n}(\mathbb{R})$ and ${\cal M}_{p\times q}(\mathbb{R}).$ Then one defines the Kronecker product of $A$ and  $B$ as being the matrix $A\oplus B$  such that
$A\oplus B=\left[\begin{array}{llll}
a_{11}B&a_{12}B&\cdots&a_{1n}B\\
a_{21}B&a_{22}B&\cdots&a_{2n}B\\
\vdots&\vdots&\cdots&\vdots\\
a_{m1}B&a_{m2}B&\cdots&a_{mn}B
\end{array}
\right].
$
}
\end{definition}
\begin{definition}
{\rm
Let $A$ e $B$ be matrices of ${\cal M}_{n}(\mathbb{R}).$ Then one defines the componentwise product of $A$ e $B$ as being the matrix $A\circ B$ such that
$(A\circ B)_{ij}=a_{ij}b_{ij},\:\forall i,j\in \{1,\cdots,n\}.$
}
\end{definition}
\begin{definition}\label{def1}
Let $B\in  {\cal M}_{p}(\mathbb{R}),$ The $kth$ kronecker power  $B^{\oplus k}$is defined inductively for all positive integers $k$  by
\begin{eqnarray*}
B^{\oplus 1}&=&B,\\
B^{\oplus k}&=&B\oplus B^{\oplus (k-1)},\:\forall k\geq 2, k\in \mathbb{N}.
\end{eqnarray*}
\end{definition}
Using the definition \ref{def1}  one considers the following notation
$$
\forall m,n\in \mathbb{N},\:(EF)^{\oplus mn}=E^{\oplus m}\oplus F^{\oplus n}
.$$
\begin{remark}
One has  $EF^{11}=E\oplus F,$  $EF^{1n}=E\oplus F^{\oplus n}$ and $EF^{n1}=E^{\oplus n}\oplus F.$
\end{remark}
\begin{definition}\label{def2}
Let $B\in {\cal M}_{p}(\mathbb{R}).$ The $kth$ componentwise  power  $B^{\oplus k}$ is defined inductively for all positive integers $k$  by
\begin{eqnarray*}
B^{\circ 1}&=&B,\\
B^{\circ k}&=&B\circ B^{\circ (k-1)}, \forall k\geq 2, k\in \mathbb{N}.
\end{eqnarray*}
\end{definition}
Using the definition \ref{def2}  one considers that
$$
\forall m,n\in \mathbb{N},\:(EF)^{\circ mn}=E^{\circ m}\circ F^{\circ n}.
$$
Let $k$ and $l$ be natural numbers.
From now on, one will use the notation that one presents next.
$$
\begin{array}{l}
\forall j\in \{1,\cdots,3\},\:E^{\circ\:jk}=(E_{j})^{\circ k},\\
\forall u,v\in \{1,\cdots,3\}\::u<v,\:E^{\circ \:uvkl}=(E_{u})^{\circ k}\circ (E_{v})^{\circ l},\\
\forall u,v\in \{1,\cdots,3\}\::u<v,\:E^{\circ \:(+uv)k}=(E_{u}+E_{v})^{\circ k},\\
\forall j,u,v\in \{1,\cdots,3\}\::u<v,\:E^{\circ \:j(+uv)kl}=(E_{j}(E_{u}+E_{v}))^{\circ kl}.
\end{array}.
$$
\begin{theorem}\label{kronecker11}
Let ${\cal M}_{p}(\mathbb{R})$ represent the vector space of the real matrices of order $p, E$ and $F$ be two idempotents matrices  of ${\cal M}_{p}(\mathbb{R}).$ Then
$E^{\oplus\:m}, F^{\oplus n}$ and $EF^{\oplus\:mn}$ are idempotents of ${\cal{M}}_{p^m}(\mathbb{R}),\: {\cal{M}}_{p^n}(\mathbb{R})$ and of ${\cal M}_{p^{m+n}}(\mathbb{R})$ respectively.
\end{theorem}
\begin{demo}
One uses induction to prove that $E^{\oplus\:m},\:F^{\oplus\:n}$ and $EF^{\oplus\:mn}$ are idempotents of ${\cal M}_{p^m}(\mathbb R),{\cal M}_{p^n}(\mathbb R)$ and of ${\cal M}_{p^{m+n}}(\mathbb R).$
\begin{enumerate}
\item[1]\label{prooff1}
$\forall m\in \mathbb{N},\:E^{\oplus\:m}$  is an idempotent of ${\cal M}_{p^m}(\mathbb R).$
\begin{itemize}
\item[i)]
Let $k=2$ and $E=[e_{ij}].$ Partitioning according to the size of $E,$
$
E^{\oplus\:2}=[e_{ij}E].
$
The $i,j$ block of $E^{\oplus\:2}E^{\oplus\:2}$ is
$$\sum^p_{k=1}e_{ik}Ee_{kj}E=\sum^p_{k=1}e_{ik}e_{kj}E=e_{ij}E.$$
But $e_{ij}E$ is the $ij$ block of $E^{\oplus\:2}.$
Then one may conclude that $$E^{\oplus\:2}E^{\oplus\:2}=E^{\oplus\:2}.$$
\item[ii)]
Let now suppose that $E^{\oplus (k-1)}E^{\oplus (k-1)}=E^{\oplus (k-1)}.$
Partitioning according to the size of  $E^{\oplus (k-1)},$  $E^{\oplus k}=[e_{ij}E^{\oplus (k-1)}].$
The $ij$ block of $E^{\oplus k}E^{\oplus k}$ is
$$\begin{array}{l}
\sum^{p}_{l=1}e_{il}E^{\oplus (k-1)}e_{lj}E^{\oplus (k-1)}\\
=\sum^{p}_{l=1}e_{il}e_{lj}E^{\oplus (k-1)}E^{\oplus (k-1)}\\
=e_{ij}E^{\oplus (k-1)}«E^{\oplus (k-1)}\\
=e_{ij}E^{\oplus (k-1)}.
\end{array}
$$
But $e_{ij}E^{\oplus (k-1)}$ is the $ij$ block of $E^{\oplus k},$
therefore $E^{\oplus k}E^{\oplus k}=E^{\oplus k}.$
\end{itemize}
One proves that $\forall n\in \mathbb{N},\;F^{\oplus\:n}$  is an idempotent matrix of ${\cal M}_{p^n}(\mathbb R)$ in the same way.
\item[2)]
$\forall m,n\in\mathbb{N},\:EF^{\oplus mn}$ is an idempotent matrix of ${\cal M}_{p^{m+n}}(\mathbb R).$

Since $EF^{\oplus mn}=E^{\oplus (m-1)}\oplus (EF^{\oplus 1n})$ one will proceed in the following way:
first one proves that $\forall n \in \mathbb{N},\: EF^{\oplus 1n}$ is an idempotent matrix of ${\cal M}_{p^{n+1}}(\mathbb{R})$ and next one proves that
$\forall m\in \mathbb{N},\: EF^{\oplus mn}$ is an idempotent matrix of ${\cal M}_{p^{m+n}}(\mathbb{R}).$
\begin{itemize}
\item[i)]$\forall n \in \mathbb{N},\: EF^{\oplus 1n}$ is an idempotent matrix of ${\cal M}_{p^{n+1}}(\mathbb{R}).$
First one proves that $$EF^{\oplus\:11}EF^{\oplus\:11}=EF^{\oplus\:11}$$
Partitioning according to the size of $F,$  $ E\oplus F=[e_{ij}F].$
The $ij$ block of $
EF^{\oplus\:11}EF^{\oplus\:11}$ is
$$\begin{array}{l}
\sum^{p}_{k=1}e_{ik}Fe_{kj}F\\
=\sum^{p}_{k=1}e_{ik}e_{kj}FF\\
=\sum^{p}_{k=1}e_{ik}e_{kj}F\\
=e_{ij}F.
\end{array}
$$
But since the $ij$  block of $EF^{\oplus:11}$ is $e_{ij}F$ then
$EF^{\oplus\:11}\oplus EF^{\oplus\:11}=EF^{\oplus\:11}$.

Now let suppose that
$$EF^{\oplus\: (1(k-1))}EF^{\oplus \:(1(k-1))}=EF^{\oplus \:1(k-1)}$$
and let prove that
$EF^{\oplus\: 1k}\oplus EF^{\oplus \:1k}=EF^{\oplus \:1k}.$
Partitioning according to the size of  $F$, $EF^{\oplus 1k}=[(EF^{\oplus 1(k-1)})_{ij}F].$
Since $$EF^{\oplus 1k}=(EF^{\oplus1(k-1)})\oplus F$$
the ij block of $$EF^{\oplus 1k}EF^{\oplus 1k}$$
is
$$\begin{array}{l}
\sum^{p^k}_{l=1}(EF^{\oplus(1(k-1))})_{il}F(EF^{\oplus(1(k-1))})_{lj}F\\
=\sum^{p^k}_{l=1}(EF^{\oplus(1(k-1))})_{il}(EF^{\oplus(1(k-1))})_{lj}FF\\
=(EF^{\oplus(1(k-1))})_{ij}F.
\end{array}.
$$
But since $(EF^{\oplus(1(k-1))})_{ij}F$ is the $ij$ block of $(EF^{\oplus 1(k-1)})\oplus F,$ this is of $E F^{\oplus 1k}$ then
$EF^{\oplus 1k}EF^{\oplus 1k}=EF^{\oplus 1k}.$
\item[ii)]
$\forall m\in \mathbb{N},\:EF^{mn}$ is an idempotent matrix of ${\cal M}_{p^{m+n}}(\mathbb{R}).$
One has already proved that
$EF^{\oplus 1n}\oplus EF^{\oplus 1n}=EF^{\oplus1n}.$
Let now suppose that
$$EF^{\oplus\:(k-1)n}EF^{\oplus\:(k-1)n}=EF^{\oplus\:(k-1)\:n}.$$
and prove that
$EF^{\oplus\:kn}EF^{\oplus\:kn}=EF^{\oplus\:kn}.$
Since $EF^{kn}=E\oplus EF^{(k-1)n}$ and partitioning according to the size of  $EF^{\oplus\:(k-1)n}$, $[EF^{\oplus kn}]=[e_{ij}EF^{\oplus\:(k-1)n}]$ then
the $ij$ block of $EF^{\oplus\:kn}$
is
$$\begin{array}{l}
\sum^{p}_{l=1}e_{il}EF^{\oplus\:(k-1)n}e_{lj}EF^{\oplus\:(k-1)n}\\
=\sum^{p}_{l=1}e_{il}e_{lj}EF^{\oplus\:(k-1)n}\\
=e_{ij}EF^{\oplus\:(k-1)n}.
\end{array}
$$
Since $e_{ij}EF^{\oplus\:(k-1)n}$  is the  $ij$ block of $EF^{\oplus\:kn}$
then
$EF^{\oplus\:kn}EF^{\oplus\:kn}=EF^{\oplus\:kn}.$
\end{itemize}
\end{enumerate}
\end{demo}
\begin{lemma}\label{kronecker12}
Let $E$ and $F$ be  idempotents matrices of ${\cal M}_{p}(\mathbb R).$ Then $E^{\circ\:m},\:F^{\circ \:n},\:EF^{\circ mn}$ are principal submatrices of
$E^{\oplus m},\:F^{\oplus \:n},\:EF^{\oplus mn}$ respectively.
\end{lemma}
\section{Necessary conditions for the existence of a strongly regular graph}\label{clara1}
Let k,l be natural numbers.
Since
$$\forall i=1,\cdots,3,\:\forall u,v,j\in \{1,\cdots,3\}\::u<v,
E^{\circ jk},
E^{\circ \:uvkl},
E^{\circ \:(+uv)k},
E^{\circ \:j(+uv)kl}
$$ are elements of the Euclidean Jordan algebra ${\cal V}_{n}$
then $\forall j=1,\cdots,3$ and $\forall u,v\in \{1,\cdots,3\}$ such that $u<v$ there exists real numbers  $q^i_{jjk},q^i_{uvkl},q^{i}_{(+uv)k}$ and $q^i_{j(+uv)kl}$  for $i=1,\cdots,3$ such that
$$
\begin{array}{l}
E^{\circ\:jk}=\sum^{3}_{i=1}q^i_{jjk}E_{i},\\
E^{\circ \:uvkl}=\sum^{3}_{i=1}q^i_{uvkl}E_{i},\\
E^{\circ \:(+uv)k}=\sum^{3}_{i=1}q^i_{(+uv)k}E_{i},\\
E^{\circ\:(j(+uv)kl)}=\sum^{3}_{i=1}q^i_{j(+uv)kl}E_{i}.
\end{array}
$$
By theorems \ref{interlace},\ref{kronecker11} and lemma \ref{kronecker12} one may conclude that,  for $i=1,\cdots,3$ and $k,l\in \mathbb{N},$ one has
$$
\begin{array}{l}
0\leq  q^{i}_{jj}\leq 1,\forall j=1,\cdots,3,\\
0\leq q^{i}_{uvkl}\leq 1,\forall u,v\in \{1,\cdots,3\}:u<v,\\
0\leq q^{i}_{(+uv)k}\leq 1,\:\forall u,v\in \{1,\cdots,3\}:u<v,\\
0\leq q^i_{j(+uv)kl}\leq 1,\forall j=1,\cdots,3,\forall u,v\in \{1,\cdots,3\}:u<v.
\end{array}
$$
After some calculations and working in the Euclidean Jordan algebra ${\cal V}_{n}$ with the componentwise product of matrices one obtains:
\begin{eqnarray*}
q^1_{11k}&=&\frac{1}{n^k}+\frac{1}{n^k}p+\frac{1}{n^k}(n-p-1),\label{desig0}\\
q^2_{11k}&=&\frac{1}{n^k}+\frac{1}{n^k}r+\frac{1}{n^k}(-r-1),\label{desig00}\\
q^3_{11k}&=&\frac{1}{n^k}+\frac{1}{n^k}s+\frac{1}{n^k}(-s-1),\label{desig000}\\
q^{1}_{22k}&=&(\frac{|s|n+s-p}{n(r-s)})^k+(\frac{n+s-p}{n(r-s)})^kp+(\frac{s-p}{n(r-s)})^k(n-p-1),\label{desig1}\\
q^{2}_{22k}&=&(\frac{|s|n+s-p}{n(r-s)})^k+(\frac{n+s-p}{n(r-s)})^kr+(\frac{s-p}{n(r-s)})^k(-r-1),\label{desig2}\\
q^{3}_{22k}&=&(\frac{|s|n+s-p}{n(r-s)})^k+(\frac{n+s-p}{n(r-s)})^ks+(\frac{s-p}{n})^k(-s-1),\label{desig3}\\
q^{1}_{33k}&=&(\frac{rn+p-r}{(n(r-s)})^k+(\frac{-n+p-r}{(n(r-s)})^kp+(\frac{p-r}{n(r-s)})^k(n-p-1),
\label{desig4}\\
q^{2}_{33k}&=&(\frac{rn+p-r}{(n(r-s)})^k+(\frac{-n+p-r}{(n(r-s)})^kr+(\frac{p-r}{n(r-s)})^k(-r-1),
\label{desig5}\\
q^{3}_{33k}&=&(\frac{rn+p-r}{n(r-s)})^k+(\frac{-n+p-r}{(n(r-s)})^ks+(\frac{p-r}{n(r-s)})^k(-s-1),
\label{desig6}
\end{eqnarray*}
\begin{eqnarray*}
q^{1}_{12kl}&=&\frac{1}{n^k}(\frac{|s|n+s-p}{n(r-s)})^l+
\frac{1}{n^k}(\frac{n+s-p}{n(r-s)})^lp+\nonumber\\
&+&\frac{1}{n^k}(\frac{s-p}{n(r-s)})^l(n-p-1),\label{desig25}\\
q^{2}_{12kl}&=&\frac{1}{n^k}(\frac{|s|n+s-p}{n(r-s)})^l+\frac{1}{n^k}(\frac{n+s-p}{n(r-s)})^lr+\nonumber
\\
&+&\frac{1}{n^k}(\frac{s-p}{n(r-s)})^l(-r-1),\label{desig26}\\
q^{3}_{12kl}&=&\frac{1}{n^k}(\frac{|s|n+s-p}{n(r-s)})^l+\frac{1}{n^k}(\frac{n+s-p}{n(r-s)})^ls+\nonumber
\\
&+&\frac{1}{n^k}(\frac{s-p}{n(r-s)})^l(-s-1),\label{desig24}
\end{eqnarray*}
\begin{eqnarray*}
q^{1}_{13kl}&=&\frac{1}{n^k}(\frac{rn+p-r}{n(r-s)})^l+\frac{1}{n^k}(\frac{-n+p-r}{n(r-s)})^lp+\nonumber
\\
&+&\frac{1}{n^k}(\frac{p-r}{n(r-s)})^l(n-p-1),\label{desig22}\\
q^{2}_{13kl}&=&\frac{1}{n^k}(\frac{rn+p-r}{n(r-s)})^l+\frac{1}{n^k}(\frac{-n+p-r}{n(r-s)})^lr+\nonumber
\\
&+&\frac{1}{n^k}(\frac{p-r}{n(r-s)})^l(-r-1),\label{desig23}\\
q^{3}_{13kl}&=&\frac{1}{n^k}(\frac{rn+p-r}{n(r-s)})^l+\frac{1}{n^k}(\frac{-n+p-r}{n(r-s)})^ls+\nonumber
\\
&+&\frac{1}{n^k}(\frac{p-r}{n(r-s)})^l(-s-1),\label{desig24}
\end{eqnarray*}
and
\begin{eqnarray*}
q^{1}_{23kl}&=&\frac{(|s|n+s-p)^k(rn+p-r)^l}{(n(r-s))^{k+l}}+\frac{(n+s-p)^k(-n+p-r)^l}{(n(r-s))^{k+l}}p+\nonumber
\\
&+&\frac{(s-p)^k(p-r)^l}{(n(r-s))^{k+l}}(n-p-1),\label{desig7}\\
q^{2}_{23kl}&=&\frac{(|s|n+s-p)^k(rn+p-r)^l}{(n(r-s))^{k+l}}+\frac{(n+s-p)^k(-n+p-r)^l}{(n(r-s))^{k+l}}r+\nonumber
\\
&+&\frac{(s-p)^k(p-r)^l}{(n(r-s))^{k+l}}(-r-1),\label{desig8}\\
q^{3}_{23kl}&=&\frac{(|s|n+s-p)^k(rn+p-r)^l}{(n(r-s))^{k+l}}+\frac{(n+s-p)^k(-n+p-r)^l}{(n(r-s))^{k+l}}s+\nonumber
\\
&+&\frac{(s-p)^k(p-r)^l}{(n(r-s))^{k+l}}(-s-1),\label{desig8}
\end{eqnarray*}
\begin{eqnarray*}
q^{1}_{(+12)k}&=&(\frac{|s|n+r-p}{n(r-s)})^k+(\frac{n+r-p}{n(r-s)})^kp
 +(\frac{r-p}{n(r-s)})^k(n-p-1),\label{desig12}\\
q^{2}_{(+12)k}&=&(\frac{|s|n+r-p}{n(r-s)})^k+(\frac{n+r-p}{n(r-s)})^kr
 +(\frac{r-p}{n(r-s)})^k(-r-1),\label{desig13}\\
q^{3}_{(+12)k}&=&(\frac{|s|n+r-p}{n(r-s)})^k+(\frac{n+r-p}{n(r-s)})^ks
 +(\frac{s-p}{n(r-s)})^k(-s-1),\label{desig14}
\end{eqnarray*}
\begin{eqnarray*}
q^{1}_{(+13)k}&=&(\frac{rn+p-s}{n(r-s)})^k+(\frac{-n+p-s}{n(r-s)})^kp
 +(\frac{p-s}{n(r-s)})^k(n-p-1),\label{desig9}\\
q^{2}_{(+13)k}&=&(\frac{rn+p-s}{n(r-s)})^k+(\frac{-n+p-s}{n(r-s)})^kr
 +(\frac{p-s}{n(r-s)})^k(-r-1),\label{desig10}\\
q^{3}_{(+13)k}&=&(\frac{rn+p-s}{n(r-s)})^k+(\frac{-n+p-s}{n(r-s)})^ks
 +(\frac{p-s}{n(r-s)})^k(-s-1),\label{desig11}\\
\end{eqnarray*}
\begin{eqnarray*}
q^{1}_{(+23)k}&=&(\frac{n-1}{n})^k+(-1)^k(\frac{1}{n})^kp
 +(-1)^k(\frac{1}{n})^k(n-p-1),\label{desig15}\\
q^{2}_{(+23)k}&=&(\frac{n-1}{n})^k+(-1)^k(\frac{1}{n})^kr
 +(-1)^k(\frac{1}{n})^k(-r-1),\label{desig16}\\
q^{3}_{(+23)k}&=&(\frac{n-1}{n})^k+(-1)^k(\frac{1}{n})^ks
 +(-1)^k(\frac{1}{n})^k(-s-1),\label{desig17}
\end{eqnarray*}
\begin{eqnarray*}
q^{1}_{3(+13)kl}&=&\frac{(rn+p-r)^k(rn+p-s)^l}{(n(r-s))^{k+l}}+\frac{(-n+p-r)^k(-n+p-s)^l}{(n(r-s))^{k+l}}p+\nonumber\\
 &+&\frac{(p-r)^k(p-s)^l}{(n(r-s))^{k+l}}(n-p-1),\label{desig18}\\ q^{2}_{3(+13)kl}&=&\frac{(rn+p-r)^k(rn+p-s)^l}{(n(r-s))^{k+l}}+\frac{(-n+p-r)^k(-n+p-s)^l}{(n(r-s))^{k+l}}r+\nonumber\\
 &+&\frac{(p-r)^k(p-s)^l}{(n(r-s))^{k+l}}(-r-1),\label{desig19}\\ q^{3}_{3(+13)kl}&=&\frac{(rn+p-r)^k(rn+p-s)^l}{(n(r-s))^{k+l}}+\frac{(-n+p-r)^k(-n+p-s)^l}{(n(r-s))^{k+l}}s+\nonumber\\
 &+&\frac{(p-r)^k(p-s)^l}{(n(r-s))^{k+l}}(-s-1),\label{desig20}
\end{eqnarray*}
\begin{eqnarray*}
q^{1}_{2(+13)kl}&=&\frac{(|s|n+s-p)^k(rn+p-s)^l}{(n(r-s))^{k+l}}+\frac{(n+s-p)^k(-n+p-s)^l}{(n(r-s))^{k+l}}p+\nonumber\\
 &+&\frac{(s-p)^k(p-s)^l}{(n(r-s))^{k+l}}(n-p-1),\label{desig21}\\ q^{2}_{2(+13)kl}&=&\frac{(|s|n+s-p)^k(rn+p-s)^l}{(n(r-s))^{k+l}}+\frac{(n+s-p)^k(-n+p-s)^l}{(n(r-s))^{k+l}}r+\nonumber\\
 &+&\frac{(s-p)^k(p-s)^l}{(n(r-s))^{k+l}}(-r-1),\label{desig22}\\ q^{3}_{2(+13)kl}&=&\frac{(|s|n+s-p)^k(rn+p-s)^l}{(n(r-s))^{k+l}}+\frac{(n+s-p)^k(-n+p-s)^l}{(n(r-s))^{k+l}}s+\nonumber\\
 &+&\frac{(s-p)^k(p-s)^l}{(n(r-s))^{k+l}}(-s-1).\label{desig23}
\end{eqnarray*}
One now make some observation about this notation.
\begin{remark}
For $i=1,\cdots,3$ the parameters $q^{i}_{jj2},$ for $j=1,\cdots,3,\:q^{i}_{uv11}$ such that $u<v\wedge u,v\in\{1,\cdots,3\}$ are the
parameters of Krein of the strongly $(n,p;a.c)$ regular graph $\tau.$
Next one must observe that $\forall i=1,\cdots,3,\forall j=1,\cdots,3,\:\:q^{i}_{jj1}=\delta_{ij}.$
From now on, for $i=1,\cdots,3$ one will call the parameters $q^{i}_{jjk},q^{i}_{(+uv)k}$ with $k\geq 3,$ and the parameters  $q^{i}_{uvkl},q^{i}_{j(+uv)kl}$ with  $l+k\geq 3$  the generalized Krein parameters of the strongly $(n,p;a,c)$ regular graph $\tau.$
\end{remark}
\begin{remark}
One now presents some consequences of the  parameters $q^{i}_{(+uv)k}$ such that $u<v\wedge u,v\in \{1,\cdots,3\}.$ Let $u$ and $v$ be natural numbers such that $u<v\wedge u,v\in \{1,\cdots,3\}.$ Since
$$
(E_{u}+E_{v})\circ (E_{u}+E_{v})\\
=E_{u}\circ E_{u}+2E_{u}\circ E_{v}+E_{v}\circ E_{v}
$$
then for $i=1,\cdots,3,$
$
q^{i}_{(+uv)2}=q^{i}_{uu}+2q^{i}_{uv}+q^{i}_{vv}.
$
But for $i=1,\cdots,3,q^{i}_{(+uv)2}\leq 1.$ Therefore one may conclude that
$$\forall i=1,\cdots,3,\:0\leq q^{i}_{uu}+2q^{i}_{uv}+q^{i}_{vv}\leq 1.$$
\end{remark}
\begin{remark}\label{observacao}
Since the generalized Krein parameters of a strongly regular graph $\tau$ are greater than zero its natural to establish necessary conditions for the existence of a strongly regular graph with  these parameters. One will analyze in this work only the generalized Krein parameters of $\tau$ associated to the eigenvalue $p$ of $\tau.$
Let $k\in 2\mathbb{N}+1$ .The parameters $q^1_{33k},q^{1}_{(+13)k}$,$q^1_{2(+13)uv}$ with $v\in 2\mathbb{N}+1$ and $q^1_{3(+13)uv}$ with $u+v\in 2\mathbb{N}+1$ permits us to establish easily criterions to see that doesn't exist a strongly $(n,p;a,c)$ regular graph.
 Each  expression of each generalized parameter  interpreted as polynomial in $n$ give us the information that if $\tau$ is a strongly $(n,p;a,c)$ regular graph such that coefficient of  the power of $n$ with exponent equal to the degree of this polynomial  is negative then one can say that if $n$ is sufficient large then this parameter is negative. Note that
 $$
 \begin{array}{l} (n(r-s))^kq^{1}_{33k}=(r^{k}-p)n^k+\sum^{k-1}_{i=0}\alpha_{i}n^{i}(k\in 2\mathbb{N}+1),\\
((n(r-s))^k)q^{1}_{(+13)k}=(r^{k}-p)n^k+\sum^{k-1}_{i=0}\beta_{i}n^{i}(k\in 2\mathbb{N}+1),\\
(n(r-s))^{u+v}q^1_{2(+13)uv}=(|s|^ur^v-p)n^{u+v}+\sum^{u+v-1}_{i=0}\delta_{i}n^{i}(v\in 2\mathbb{N}+1),\\
(n(r-s))^{u+v}q^1_{3(+13)uv}=(r^{u+v}-p)n^{u+v}+\sum^{u+v-1}_{i=0}\gamma_{i}n^{i}(u+v\in 2\mathbb{N}+1).
\end{array}
$$
\end{remark}
The remark \ref{observacao} conduct us to the theorem \ref{kronecker16}.

\begin{theorem}\label{kronecker16}
Let  $\tau$ be a strongly $(n,p;a,c)$ regular graph such that $0<c<p<n-1$ with adjacency matrix $A$, having eigenvalues $p,r$ and $s.$ Then
\begin{eqnarray*}
(rn+p-r)^k&+&(-n+p-r)^kp+(p-r)^k(n-p-1)\geq 0\nonumber\\
\forall k\in &2&\mathbb{N}+1\:\quad(q^{1}_{33k}\geq 0),
\label{desig30}\\
(rn+p-s)^k&+&(-n+p-s)^kp+(p-s)^k(n-p-1)\geq 0\nonumber\\
\forall k\in &2&\mathbb{N}+1\label{desig31}\:\quad (q^{1}_{(+13)k}\geq 0),\nonumber\\
\\
(rn+p-r)^k&\!&(rn+p-s)^l+(-n+p-r)^k(-n+p-s)^lp+\nonumber\\
+(p-r)^k&\!&(p-s)^l(n-p-1)\geq 0 \nonumber \\
 \forall k,l\in \mathbb{N},l+k &\in&2\mathbb{N}+1\quad (q^1_{3(+13)kl}\geq 0)\label{desig32},\\
(|s|n+s-p)^k&\!&(rn+p-s)^l+(n+s-p)^k(-n+p-s)^lp+\nonumber\\
 +(s-p)^k&\!&(p-s)^l(n-p-1)\geq 0\nonumber \\
 \forall k\in \mathbb{N},l &\in& \quad 2\mathbb{N}+1 (q^1_{2(+13)kl}\geq 0).\label{desig33}
\end{eqnarray*}
\end{theorem}
One presents the lemma  \ref{kronecker17}, since the left member of each inequality\footnote{For instance, from the first inequality of lemma \ref{kronecker17} one may establish the corollary \ref{kronecker15}}  of this lemma is a polynomial in  $\mathbf{n}$ of degree three is more suitable for getting information, like that presented in corollary \ref{kronecker15}, about a strongly regular graph.
\begin{lemma}\label{kronecker17}
Let $\tau$ be a strongly $(n,p;a,c)$ regular graph such that $0<c<p<n-1$  with adjacency matrix $A$, having eigenvalues $p,r$ and $s.$ Then
\begin{eqnarray}
(rn+p-r)^3&+&(-n+p-r)^3p+(p-r)^3(n-p-1)\geq 0\nonumber\\
\quad &\!&(\quad q^{1}_{333}\geq 0),
\label{desig34}\\
(rn+p-s)^3&+&(-n+p-s)^3p+(p-s)^3*(n-p-1)\geq 0\nonumber\\
\quad &\!&(q^{1}_{(+13)3}\geq 0),\nonumber\\
\\
(rn+p-r)^2&\cdot&(rn+p-s)+(-n+p-r)^2(-n+p-s)p+\nonumber\\
(p-r)^2&\cdot&(p-s)(n-p-1)\geq 0, \nonumber \\
 \quad &\!&(q^1_{3(+13)21}\geq 0),\label{desig35}\\
 (rn+p-r)&\cdot&(rn+p-s)^2+(-n+p-r)(-n+p-s)^2p+\nonumber\\
 (p-r)&\cdot&(p-s)^2(n-p-1)\geq 0\:\quad (q^{1}_{1312}\geq 0)\nonumber\\
 \quad &\!&(q^1_{3(+13)12}\geq 0),\label{desig36}\\
(|s|n+s-p)^2&\cdot&(rn+p-s)+(n+s-p)^2(-n+p-s)p+\nonumber\\
 +(s-p)^2&\cdot&(p-s)(n-p-1)\geq 0\nonumber \\
 \quad &\!&(q^1_{2(+13)21}\geq 0).\label{desig37}
\end{eqnarray}
\end{lemma}
From lemma \ref {kronecker17} and (\ref{desig34}) one  obtains the  corollary \ref{kronecker15}.
\begin{corollary}\label{kronecker15}
If $\tau$ is a strongly $(n,p;a,c)$ regular graph and $r$ is the positive eigenvalue of his matrix of adjacency then:
\begin{enumerate}
\item[$1$]
$$r<  p^{\frac{1}{3}}$$
$$\Downarrow$$ $$n<\frac{(p-r)(3r^2+3p+\sqrt{r^4+18pr^2+p^2+8r^3p+8pr})}{2(p-r^3)}$$
\item[$2$]
$$r> p^{\frac{1}{3}}$$
$$\Downarrow$$
$$n>\frac{(p-r)(3r^2+3p+(p-r)\sqrt{r^4+18pr^2+p^2+8r^3p+8pr})}{2(p-r^3)}$$
\end{enumerate}
\end{corollary}.

\begin{conclusion}
Each necessary condition for the existence of a strongly  $(n,p;a,c)$ regular graph of theorem \ref{kronecker16} allows us to conclude\footnote{For instance, interpreting (\ref{desig30}) as a polynomial in $n$ one has $$(rn+p-r)^k+(-n+p-r)^kp+(p-r)^k(n-p-1)=(r^{k}-p)n^k+\sum^{k-1}_{i=0}\alpha_{i}n^{i}$$ and so is natural to conclude that, when k is an odd natural number, if $r^k-p<0$ and $n$ is big then $(rn+p-r)^k+(-n+p-r)^kp+(p-r)^k(n-p-1)<0.$} that $n$ can not be too big when the coefficient of the power of $n$ with greatest exponent is negative.
One presents the lemma \ref{kronecker17}, since the left member of each inequality of this lemma \footnote{For instance, from the first inequality of lemma \ref{kronecker17} one may establish the corollary \ref{kronecker15}}   is a polynomial in  $\mathbf{n}$ of degree three and therefore is more suitable for getting information, like that presented in corollary \ref{kronecker15}, about a strongly regular graph.
Finally one observes  that the introduction of the generalized  Krein  parameters of a strongly  $(n,p;a,c)$ regular graph allowed us to conclude that the Krein parameters for $i,u,v\in \{1,\cdots,3\}:u<v$ must satisfy not only $0\leq q^{i}_{uv}\leq 1\wedge 0\leq q^{i}_{vv}\leq 1$ but also $q^{i}_{uu}+2q^{i}_{uv}+q^{i}_{vv}\leq  1.$
\end{conclusion}

\end{document}